\documentstyle[11pt]{article}
\title{{\bf Comment on a Paper by Kuo, Parusinski and Paunescu On Jacobian Conjecture}}
\author{{\Large T. Moh}\thanks{ Math Department, Purdue University, West Lafayette, Indiana 47907-1395. tel: (765)-494-1930, e-mail ttm@math.purdue.edu}}
\date{Oct 10, 2005}
\begin{document}
\maketitle

\begin{abstract}
The said paper entitled "A Proof Of The Plane Jacobian Conjecture" is not true.
\end{abstract}

\noindent{\bf Comments}

We shall use the following notations: Let us call the paper by Kuo, Parusinski $\&$ Paunescu
entitled {\bf A PROOF OF THE JACOBIAN CONJECTURE} Sept 2005, http://math.univ-angers.fr/$\sim$parus/
as [Kuo], and the paper by Moh entitled {\bf On the Jacobian Conjecture} Crelle 1983
pp 140-212 as [Moh].

J.Lipman sent us the URL of [Kuo]. Since we are busy in revising our lecture note on
{\it Algebraic Coding Theory}, we ask A. Sathaye, a very good mathematician, for opinion.  
In an e-mail sent by  A. Sathaye to us, he states,

{\it I have looked at the proof without going thru all the details.
They seem to claim that there are no "minor roots" in your language.
(That is what their "eclipse root" is!)$\cdots\cdots$, I cannot point to a hole,
but I don't have much confidence in the proof at this point.}

This short note is an explanation of A. Sathaye's words. T.C.Kuo, a known mathematician 
in analytic functions theory and Jet theory,  and his collaborators apparently were unaware 
of our paper in 1983. 

There are many similarities between [Kuo] and [Moh]. For instance, we have the
following partial table,

\begin{center}
\begin{tabular}{|r|r|r|r|}\hline
[Kuo]&ref. &[Moh] & ref. \\[0.5ex]
\hline\hline
$\gamma^{*}$&p.9&$\pi$-root&p.146 \\
\hline
proportional&p.2&non-splitting&p.163\\
\hline
resonant&p.4&splitting&p.163\\
\hline
weak resonant&p.4&---&p.206\\
\hline
\end{tabular}
\end{center}

We may point most Propositions of [Kuo]  to the similar ones in [Moh], and
some of Propositions of [Kuo] are unique (hence maybe wrong). Let
us examine one of them, Proposition 4.4 (Spider web lemma) on p.10 of
[Kuo]. In the proof, the authors claim (on the 11th line of the proof) 
that :{\it Hence J(F,G)$\equiv$ 0, ($p,\bar{p}$) and ($q, \bar{q}$)
are proportional.} First, the statement is meaningless if given two
arbitrary polynomials $f(z)$ and $g(z)$ and 
\begin{eqnarray*}
&&F(x,y)=f(x^ay^{-b})+\cdots\\
&&G(x,y)=g(x^ay^{-b})+\cdots
\end{eqnarray*}
Then the initial form of $J(F,G)\equiv 0$. Certainly one can not deduce that
their Newton polygons $N_1, N_2$ are {\it radially proportional} (p.10).

This kind mistake was first observed by S.S.Abhyankar in the early 1970's
about a claim of Jacobian Conjecture by Mr. Jan. Namely, one may solve
the Jacobian equation at $\infty$ to produce a contradiction. However,
the trouble is after a few steps of solving, the two power series become
two units (with order $0$'s), say

\begin{eqnarray*}
&&f(xy^a+\gamma,y^{-1})=f_0(x)+f_{\delta}(x)y^{\delta}+\cdots\\
&&g(xy^a+\gamma,y^{-1})=g_0(x)+g_{\epsilon}(x)y^{\epsilon}+\cdots
\end{eqnarray*}
Then the Jacobian equation is quite different from the previous cases.
An experienced mathematician may easily find the similarity of the above
arguments and our argument in the note which substantiate the point of view
of A. Sathaye.

As for Kuo and his collaborators, we believe that they have a good taste
of mathematics, and wish that they  will push the analytic method deeper to 
solve the Jacobian Conjecture.

\vskip 0.1 in

\noindent{\bf Remark (on Dec 22, 2005)}

This note was written on Oct 10, 2005 and was sent to the authors. At once they replied to insist that
they are correct, which was natural. After a month we checked the website of Parusinski, and found that
a new sentence "{\it The proof contains some gaps in section 7"} by the authors without mentioning any
objection by us.

\end{document}